\documentclass[%
 aip,
rsi,%
 amsmath,amssymb,
 reprint,%
]{revtex4-1}

\usepackage{graphicx}
\usepackage{dcolumn}
\usepackage{bm}
\usepackage{xcolor}

\graphicspath{{Figures/}}

\newcommand{\iext}{I}

\newcommand{\tstim}{T_{\rm stim}}

\newcolumntype{P}[1]{>{\centering\arraybackslash}p{#1}}

\begin{document}

\preprint{AIP/123-QED}

\title{Order-indeterminant event-based maps for learning a beat}

\author{\'{A}. Byrne}
 \email{aine.byrne@ucd.ie}
 \affiliation{School of Mathematics and Statistics, University College Dublin, Ireland}
 \affiliation{Center for Neural Science, New York University, New York, NY 10003 USA}
\author{J. Rinzel}
 \affiliation{Center for Neural Science, New York University, New York, NY 10003 USA}
 \affiliation{Courant Institute of Mathematics, New York University, New York, NY 10003 USA}
\author{A. Bose}%
\affiliation{ 
Department of Mathematical Sciences, New Jersey Institute of Technology, Newark, NJ 07102 USA
}%

\date{\today}

\begin{abstract}
The process by which humans synchronize to a musical beat is believed to occur through error-correction where an individual's estimates of the period and phase of the beat time are iteratively adjusted to align with an external stimuli.  Mathematically, error-correction can be described using a two-dimensional map where convergence to a fixed point corresponds to synchronizing to the beat. In this paper, we show how a neural system, called a \emph{beat generator}, learns to adapt its oscillatory behaviour through error-correction to synchronize to an external periodic signal. We construct a two-dimensional event-based map which iteratively adjusts an internal parameter of the beat generator to speed up or slow down its oscillatory behaviour to bring it into synchrony with the periodic stimulus. The map is novel in that the order of events defining the map are not {\it a priori} known.  Instead, the type of error-correction adjustment made at each iterate of the map is determined by a sequence of expected events. The map possesses a rich repertoire of dynamics, including periodic solutions and chaotic orbits. 
\end{abstract}

\maketitle

\begin{quotation}
Music is an important part of human society.
As humans, we have an innate ability to quickly recognise rhythmicity and reproduce it by moving our body to the music. To bring our movements into sync with the music, we make a series of adjustments to speed up or slow down in a process  known as error-correction. The time between successive movements, such as claps, serves as an estimate of the beat period, and is 
compared to the actual time between musical beats. The exact timing of each clap compared to the beat time determines whether or not our clapping phase is aligned with the beat. Using these two comparisons we make a judgement on whether we need to speed up or slow down our clapping. 
In this paper, we explore how such an error-correction scheme could be implemented by a neuronal network. Our work focuses on the derivation and analysis of a two-dimensional map to describe how an error-correction scheme can bring a neural system into sync with a periodic external drive. Importantly, the map can describe situations in which the order of events changes. This is critical because we often aperiodically alternate between being too early and too late when learning a beat, i.e. on one cycle we may clap just before the beat, while on the next we clap after it. We call this property order-indeterminacy.
\end{quotation}

\section{\label{sec:intro} Introduction}

Humans can easily recognize rhythmicity within speech and music which spans the range of 0.5 to 10 Hz \cite{Grahn2012}. The ability to discern and track a periodic structure in music is called {\it beat perception}. While a piece of music may be quite complex in terms of its rhythmic structure, experimental studies have shown that neurons in various parts of the brain exhibit oscillations in their voltage profiles that match the beat frequency \cite{Henry2017,Merchant2015}.  A primary way to assess beat perception is through finger tapping experiments where participants tap at what they perceive to be the beat or along to an isochronous (evenly spaced in time) metronome \cite{Repp2005,Repp2013}. Intertap intervals are compared to interbeat intervals to measure the participants ability to match the period. The exact timing of taps at each cycle is compared to beat times to measure the phase of tapping. These two measures, period and phase, have led to a set of error-correction models \cite{Mates1994,Mates1994a}. These iterative algorithms attempt to describe how humans use error measurements at each cycle to correct their tap times in subsequent cycles. In essence, the schemes defines a two-dimensional, event-based map.

Event-based maps refer to a class of dynamical systems where the set of dependent variables are updated at the occurrence of a particular event, such as a trajectory crossing a Poincar\'e section. Event-based maps arise in other contexts besides error-correction. In neural systems, maps based on spike timing are often derived. For example, using a phase response curve in the presence of a weak-coupling assumption between neurons is used to assess the existence and stability of phase-locked solutions \cite{Achuthan2009, Oprisan2010}. In cardiac systems, maps have been used to study the repolarization of the left and right ventricles\cite{Strumillo2002}.  In the field of robotic movement, event-based maps are used to implement specific control strategies for stabilization of walking \cite{Hamed2013}.

In this work, we derive a map that originates from our recently developed mathematical model for how humans learn to generate a beat \cite{Bose2019}. We proposed a biophysical framework whereby a neuronal oscillator, called a {\it beat generator} (BG), learns the period and phase of an external, isochronous tone sequence (S). We derived an error-correction scheme that adjusts a biophysical parameter on a cycle-by-cycle basis allowing the BG to synchronize its voltage spikes to the tone times. Updates to this parameter occur at BG spikes and stimulus tones, as such they are characterized as events, and form the basis of the event-based map. Our previous results \cite{Bose2019} are largely numerical, but the simulation results give rise to set of mathematical questions that we shall address here. 

In the ideal case, the learning process can be described in a 1-to-1 manner in which every BG spike is followed by an S tone and vice versa. This leads to a situation of monotone convergence to the synchronized solution in which the order of events is preserved. However, when tapping to a beat, our tap times typically jitter in a neighborhood of the actual stimulus tone onset times. For example, if at one cycle, the tap is too late, we might compensate by tapping at a faster rate. This  may then result in the second tap occurring before the second tone. These situations can be described as order-indeterminant because at the moment of any one event, say a tap, it is not {\it a priori} known whether the next event will be a tone or a tap. The event itself depends on the adjustments made on the fly by the listener in an attempt to synchronize to the tone sequence.  Modeling this situation mathematically is a challenge as one can no-longer make a 1-to-1 (order preserving) assumption about events. Instead, we derive a novel two-dimensional order-indeterminant event-based (OIEB) map that can predict whether and how the BG will synchronize with the stimulus, even when the order of events is not prescribed from cycle-to-cycle. The main difficulty in deriving this map is that different adjustments are made at stimulus tones and BG spikes. Thus when the order it unknown, the sequence of updates is unknown. One of the important contributions of this paper is to show how to systematically overcome this obstacle to make the appropriate number and type of updates per cycle.

An interesting property of the OIEB map is that it is piecewise-smooth, with a discontinuity arising due to the update rule for phase learning. As such, our study falls within the category of low-dimensional discontinuous maps, which are seen in a variety of contexts. From a theoretical point of view, some of the earliest studies of this type involved piecewise linear one-dimensional maps \cite{Leonov1962}. Discontinuous maps arise more generally in switching systems, where the vector field changes discontinuously across a lower dimensional manifold \cite{DiBernardo2008}. Biological applications have provided an abundance of examples. In particular, systems defined on periodic phase spaces such as a circle or torus lead to discontinuous maps. Applications of such maps include cardiac dynamics \cite{Glass1991}, circadian rhythms \cite{Diekman2016}, and sleep rhythms \cite{Booth2017}, to name just a few. 

The organization of this paper is based on deriving the OIEB map from a series of intermediate steps to build maps of increasing complexity. The simplest such map is one-dimensional and accounts just for period learning. We will show the conditions under which this one-dimensional map has a stable fixed point which corresponds to a learned period and then describe how the loss of stability leads to higher periodicity orbits. The next level of complexity involves phase learning, leading to the two-dimensional order preserving map. Here, we conduct a linear stability analysis for the synchronized (learned) solution. We demonstrate that while there are  regions of parameter space in which the synchronized solution is stable, only part of the parameter space corresponds to order preserving solutions. This leads us to derive the OIEB map where information from the stability analysis in the order preserving case will be used to explain the observed dynamics of the OIEB map. This map possesses a rich set of dynamics that include not just the existence of fixed points, but also various periodic orbits, as well as chaotic solutions. 


\section{Model for Beat Generation}
We briefly introduce the important aspects of our beat generation model \cite{Bose2019}. We first define a periodic stimulus tone sequence with an initial tone time, $t_0$ and interstimulus onset interval of $\tstim$, such that tones occur at $t=t_0+k\tstim$, where $k=0,1,2\ldots$. We then define a beat generator (BG) as a neuronal oscillator with easily identifiable spike times whose interspike interval can be controlled by parameters. There are a variety of biophysical models based on the Hodgkin-Huxley formalism that we could choose. Here, we utilize perhaps the simplest, the leaky integrate-and-fire ($LIF$) model, to allow us to focus on the mathematical results. The $LIF$ is given by 
\begin{equation}
\label{eq:LIF}
v'= (\iext-v)/\tau,
\end{equation}
with the reset condition $v\rightarrow0$, when $v$ crosses the firing threshold, set (without loss of generality) at 1. At the moment of reset, the BG is said to exhibit a spike. The variable $v$ represents voltage, the parameter $I$ provides an external drive, and $\tau$ is a time constant. For $I <1$, the system has a stable fixed point at $v=I$. The discontinuous reset condition leads to a bifurcation at $\iext = 1$, and the system exhibits oscillatory behaviour for $\iext > 1$. The period of these oscillations is given by
\begin{equation}
\label{eq:T}
T(I)=\tau\ln \left(\frac{\iext}{\iext-1}\right),
\end{equation} 
which is a one-to-one invertible function on the domain $I>1$, with inverse
\begin{equation}
\label{eq:I}
I(T)=\frac{1}{1-e^{-T/\tau}}.
\end{equation}
Thus each value of the period results from a unique value of the drive and vice versa.
In the map, the value of $I$ is updated at each BG spike and at each stimulus tone time, until the period and spike times of the BG exactly match those of the stimulus.
Hence, $I$ is a variable of the two-dimensional map. 
To update $I$ we must keep track of the time between consecutive BG spikes (period) and the time between the BG spike and stimulus tone (phase). The phase of the BG $\phi$ is the second variable of the map. Phase is updated at every tone time. Thus, the two-dimensional map will exhibit a type of asynchronous updating where one of the variables, $I$, is updated more frequently, than the other, $\phi$. 

\begin{figure*}
\includegraphics[width=1\textwidth]{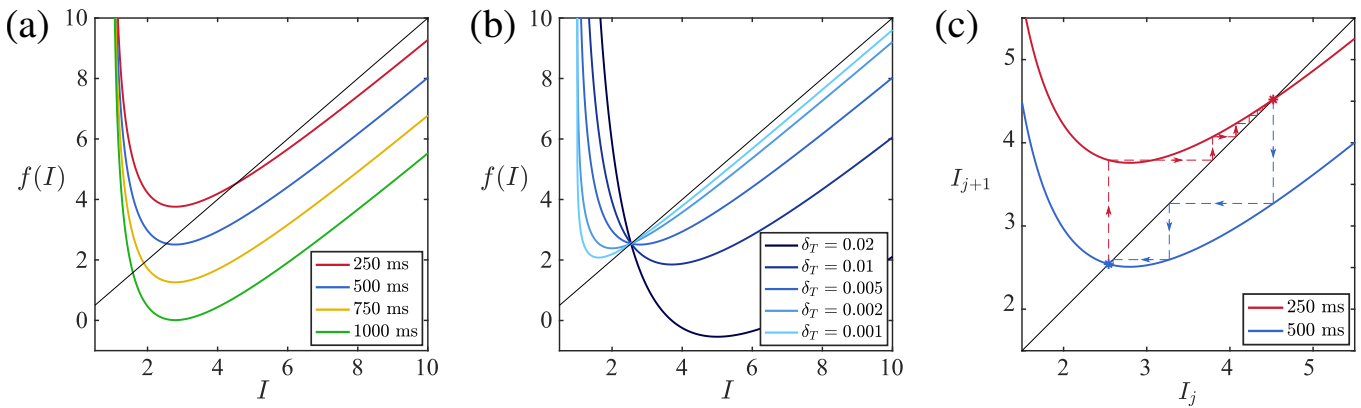}
\caption{\label{fig:period-map} \textbf{Period learning is stable provided learning rates are sufficiently small.} (a) The 1D period correction map for different stimulus periods with $\delta_T=0.005$. Each curve crosses the diagonal
at exactly one point, corresponding to a fixed point. As the stimulus period is reduced, the fixed point moves up the diagonal. Additionally, the slope at the fixed point increases as the period is reduced, becoming unstable for $\tstim=1000$ ms. (b) The 1D map for a fixed period ($\tstim=500$ ms) and varying the learning rate $\delta_T$. Increasing $\delta_T$ increases the slope at the fixed point but does not change the location of the fixed point. Stability is lost when the slope becomes less than -1 ($\delta_T=0.02$). (c) Cobweb diagram showing the iterates for a transition from $\tstim=500$ ms to $\tstim=250$ ms (red) and $\tstim=250$ ms to $\tstim=500$ ms (blue). Convergence to $\tstim=500$ ms requires significantly less iterates as the slope at the fixed point is close to zero.}
\end{figure*}

\section{\label{section-exactmap}Order Preserving 2-dimensional map}
We begin by deriving the order preserving map under the assumption of 1-to-1 firing. Namely, every spike of the $BG$ is followed by a stimulus tone and vice versa. The goal is to show the process by which the $BG$ learns the exact spike times of the stimulus as well as the interonset interval between tones. To do so, we introduce two learning rules, a period correction rule and a phase correction rule. Both rules update the value of $I$, which we treat as a variables of the map. Period correction updates are done at every spike of the BG, while phase correction updates are performed at every stimulus tone time. In Section \ref{section-oiebmap}, we will extend these results to derive the OIEB map. 

\subsection{\label{period-map}Period correction leads to a one-dimensional map}
Consider a periodic stimulus oscillating with a period $T_{stim}$ with initial spike at $t_0=0$. For the BG to match the stimulus period, it must learn the value of $\iext$ associated with that period. We define the following one-dimensional, period correction map,
\begin{equation}
\label{eq:period-map1}
I_{j+1} =  I_j + \delta_T\left[T(I)-T_{stim}\right],
\end{equation}
which iteratively increases (decreases) $\iext$ to decrease (increase) the oscillatory period of the $BG$, until it is the same as the stimulus period $T_{stim}$. The parameter $\delta_T$ is the strength of the period correction rule. Updates are associated with the spike times of the $BG$. Namely, once the update to $I_{j+1}$ is made, the next interspike interval is uniquely determined by this new value.  Using \eqref{eq:T} we rewrite \eqref{eq:period-map1} as
\begin{equation}
\label{eq:period-map2}
f(I) = I + \delta_T\left[\tau\ln \left(\frac{I}{I-1}\right)-\tstim\right],
\end{equation}
and $ I_{j+1}\equiv f(I_j)$.

The map has a unique fixed point at $I^*=1/(1-\exp(-\tstim/\tau))$, whose stability is determined using $f^\prime(I^*) = 1-\delta_T\tau/(I^*(I^*-1))$. The fixed point is stable for $|f^\prime(I^*)|<1$, and, as such, $\delta_T<2I^*(I^*-1)/\tau$. Given that $T$ is a monotonically decreasing function of $I$, it follows that to converge to a stable fixed point at larger periods, smaller $\delta_T$ is required.
The graph of $f(I)$ has a vertical asymptote at $I=1$, a local min at $I=(1+\sqrt{1+4\tau\delta_T})/2$. It is concave up and in the limit as $I\to\infty$ has a slant asymptote at $I-\delta_T T$. The slope of $f(I)$ is less than 1 for finite $I$ and tends to 1 as $I\to\infty$. Thus the graph intersects the diagonal at exactly one point, corresponding to the unique fixed point. As $\tstim$ increases, the map shifts down in the $I$-$f(I)$ plane (Fig. \ref{fig:period-map}(a)), and the unique fixed point becomes unstable when the slope at the intersection decreases below -1. For a fixed period, increasing $\delta_T$ decreases the slope of the curve at the fixed point (Fig. \ref{fig:period-map}(b)). Combining these two findings, one notes that there are privileged parameter pairs that lead to fast convergence. For example, by choosing a $\delta_T=0.005$, the slope at the fixed point associated with the $\tstim=500$ ms map is zero. For this same value of $\delta_T$, the slope of the $\tstim=250$ ms is positive, but less than 1. Thus a tempo change from $\tstim=250$ to $500$ ms will converge much faster than a change from 500 to 250 ms (Fig. \ref{fig:period-map}(c)).

\begin{figure}
\includegraphics[width=0.45\textwidth]{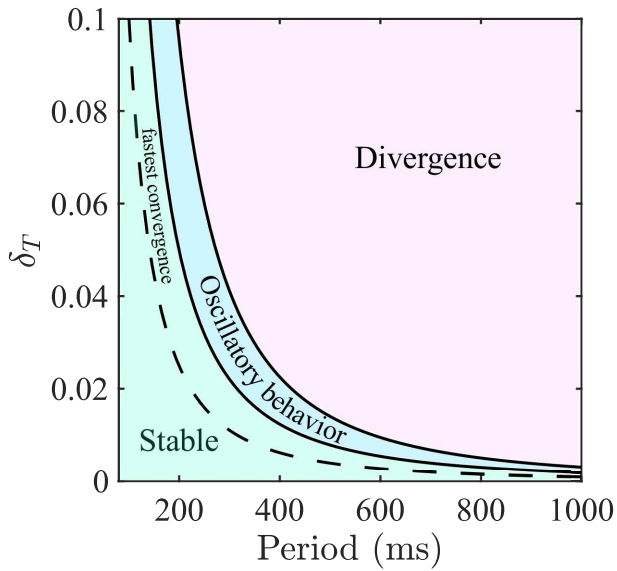}
\caption{\label{fig:delta-period} \textbf{Critical $\delta_T$ values as a function of period}.
The map converges to the to the synchronous solution ($T(I)=\tstim$) for parameter values under the first solid curve (green) and diverges for parameter values above the second solid curve (pink). Between these two regions we see oscillatory solutions, including chaotic attractors (blue). Convergence is fastest when the slope at the fixed point is zero (dashed curve).}
\end{figure}

The fixed point remains stable so long as $f^\prime(I^*)>-1$. The $\delta_T$ value for which stability is lost can be calculated by solving $f^\prime(I^*)=-1$, for $\delta_T$ (Fig. \ref{fig:delta-period}). Convergence to the fixed point is fastest when $f^\prime(I^*)=0$. Hence, solving $f^\prime(I^*)=0$ for $\delta_T$, gives the optimal $\delta_T$ value for a given period (dashed line). Notice that these values are highly dependent on the stimulus period, with smaller periods allowing for significantly larger $\delta_T$ values consistent with what is shown in Fig. \ref{fig:period-map}.  The fixed point loses stability when $f(I^*) < -1$ and in that case it is possible to obtain periodic and bounded oscillatory solutions of the map provided that the local minimum of $f(I)$ is greater than 1 (blue region). This necessary condition on the minimum of $f(I)$ guarantees that iterates of the map fall in the domain of the map, $(1,\infty)$. Periodic solutions are obtained in the usual manner, namely by searching for fixed points of suitable composition of the map with itself, e.g. period-2 solutions arise as solutions of $f^2(I)= f(f(I))=I$. For parameter values that lie in the blue region we were able to find period-2, 4 etc points. Moreover we found values of parameters at which the solution does not converge for any initial value, indicative of chaotic behavior. 
To illustrate the existence of chaos for this map, consider the period doubling route to chaos \cite{May1976}. We calculated  the ratios of the differences of parameters of successive period doubling bifurcation values  (period-1 to period-2, period-2 to period-4 , etc). These ratios turn out to be approximately $F_3=4.328$, $F_4= 4.619$ and $F_5=4.655$. Thus the period doubling behavior suggests that the value $F_n$ converges to the value $F\approx 4.669$,  the so called Feigenbaum ratio \cite{Feigenbaum1978}, giving strong indication that the one-dimensional period correction map does exhibit chaos. The findings presented here, will be instructive when we consider the dynamics of the OIEB map in later sections.  

\begin{figure*}
\includegraphics[width=1\textwidth]{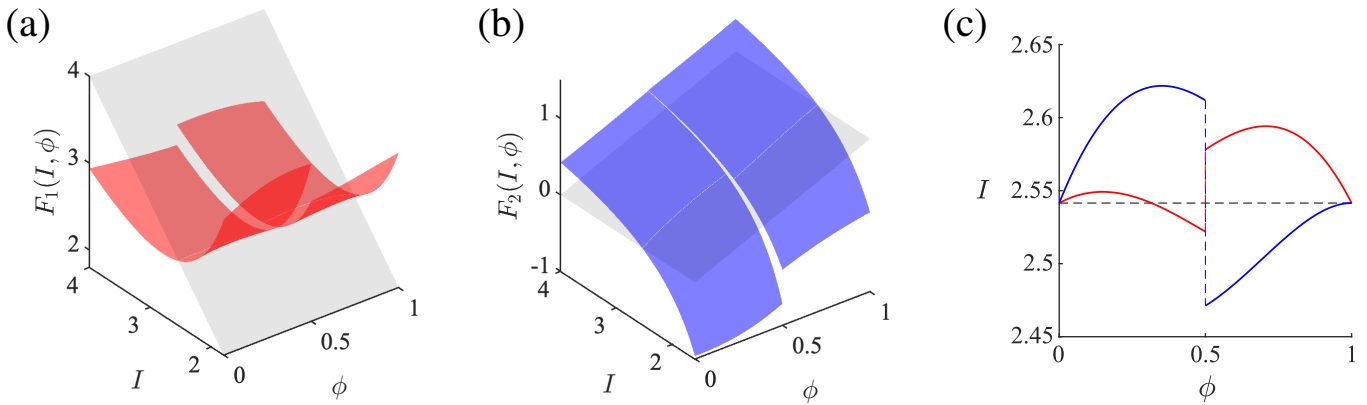}
\caption{\label{fig:period-phase-surfaces} \textbf{Fixed points of the period and phase map correspond to the synchronous solution.} (a) 
The surface $F_1(I,\phi)$ (red), which corresponds to updates of the $I$ variable in the 2D map, and the surface $z=I$ (grey). There are no updates to $I$ where the two surfaces intersect. There is a discontinuity at $\phi=0.5$ due to the sgn function in the phase update rule. 
(b)  The surface $F_2(I,\phi)$ (blue), corresponding to updates to $\phi$ in the 2D map, and the surface $z=\phi$ (grey). The phase $\phi$ is constant where the two surfaces intersect. The discontinuity at $\phi=0.5$ appears here also.
(c) A projection of the intersections for the $I$ (red) and $\phi$ (blue) surfaces with $z=I$ and $z=\phi$, respectively. Where the red and blue curves intersect, at $\phi=0$ and $\phi=1$,  correspond to the fixed points. 
Parameter values given in Appendix \ref{appendix:parameters}.}
\end{figure*}

\subsection{\label{phase-map}Phase correction leads to a 2-dimensional map}
The period rule brings the stimulus and BG periods into alignment, but pays no regard to the phase of the BG events relative to the stimulus events. As such, an additional rule is required, a phase correction rule. 
This rule also targets the drive $I$ in order to speed up/slow down the BG oscillations. It works in concert with the period correction rule, sequentially making adjustments to $I$.
We define the phase of the BG $\phi$ as the time since the last BG event divided by the stimulus period.
If $\phi<0.5$ when a stimulus event occurs, the BG is interpreted to be firing before the tone and needs to be slowed down, decreasing both  $I$ and $\phi$. If $\phi>0.5$, the opposite interpretation is taken and the BG should speed up, increasing $I$ and $\phi$. The maximal correction occurs close to 0.5. The phase update rule occurs at every stimulus event and will update the $I$ value by an amount $\delta_\phi  \Delta I_\phi(\phi)$ where 
\begin{equation}
\label{eq:DeltaIphi}
\Delta I_\phi(\phi) =   {\rm sgn}(\phi-0.5)\phi(1-\phi),  
\end{equation}
${\rm sgn}(x)$ is the sign function and $\delta_\phi>0$ is the strength of the update rule. 

The two-dimensional map requires that the $I$ value is updated at both BG spikes (period correction) and S tones (phase correction). The phase $\phi$ is only updated once per cycle, at a BG spike.
A cycle is defined as an oscillatory cycle of the BG, the $n$th cycle begins at the $n$th BG spike and ends at the $n+1$th spike. 
A complete cycle of the 2-dimensional map takes $(I_n, \phi_n)$ to $(I_{n+1}, \phi_{n+1})$ through a sequence of intermediate steps.
At the S tone, phase correction is applied
\begin{equation}
\label{eq:phase-rule}
I_{temp} = I_n + \delta_\phi \Delta I_\phi(\phi_n).
\end{equation}
The BG voltage \eqref{eq:LIF} evolves with this new value $I_{temp}$ until its next spike. 
At the BG spike, $I$ is updated again using the period rule. To apply the period rule, we must calculate the BG cycle period $T_n$.
This cycle period is comprised of two parts, the time between the $n$th BG spike and the $n$th tone, and the time between the $n$th tone and the $n+1$th BG spike,
\begin{align}
\label{eq:tnpok1}
T_n(I,\phi)=\phi\tstim +\tau \ln\left.\frac{I + \delta_\phi\Delta I_\phi(\phi)  - v(I,\phi)}{I+ \delta_\phi\Delta I_\phi(\phi) -1}\right.,
\end{align}
where
\begin{align}
\label{eq:vtnk}
v(I, \phi)=I\left(1-e^{-\tstim\phi/\tau}\right).
\end{align}
and the second term in \eqref{eq:tnpok1} is obtained by integrating (\ref{eq:LIF}) from the initial value $v(I_n, \phi_n)$ to 1.
Then the period correction is given by
\begin{align}
\label{eq:doubleperiod-rule}
I_{n+1} &= I_{temp} + \delta_T\Delta I_T(I_n,\phi_n),
\end{align}
where $\Delta I_T(I,\phi) = \left[T_n(I,\phi)-T_{stim}\right]$.
At the $(n+1)$th BG spike, we also update the phase for the next cycle,
\begin{equation}
\label{eq:phase-update}
\phi_{n+1}=\phi_n +\frac{\tstim-T_n(I_n,\phi_n) }{\tstim}.
\end{equation}

The system can then be written as a two-dimensional map, with updates to $I$ and $\phi$ on every iteration defined by the paired event of an S tone and a BG spike,
\begin{align}
\label{eq:2dmap1}
I_{n+1} &=F_1(I_{n},\phi_n)\\
\label{eq:2dmap2}
\phi_{n+1} &=F_2(I_{n},\phi_n).
\end{align}
The functions $F_1$ and $F_2$ are given by
\begin{align}
\label{eq:2d-map-i}
F_1(I,\phi) &= I+ \delta_\phi \Delta I_\phi(\phi) + \delta_T\Delta I_T (I,\phi)\\
\label{eq:2d-map-phi}
F_2(I,\phi) &= \phi + \frac{\tstim-T_n(I,\phi) }{\tstim}.
\end{align}
A fixed point of the map satisfies the algebraic conditions $I=F_1(I,\phi)$ and $\phi=F_2(I,\phi)$. 
For $\tstim=500$, we graph the surfaces $F_1(I\phi)$ and $F_2(I,\phi)$ in separate three dimensional spaces (Fig. \ref{fig:period-phase-surfaces}(a) and (b)), and examine their intersection with the two-dimensional planes $z=I$ and $z=\phi$, respectively.
The ensuing curves of intersection are then projected onto the $\phi-I$ domain (Fig. \ref{fig:period-phase-surfaces}(c)) and their intersection yields the fixed points of the map. As expected, the fixed points are located at $(I^*(\tstim),0)$ and $(I^*(\tstim),1)$. 
Though the fixed points at $\phi=0$ and $\phi=1$ correspond to the same synchronous solution, we will show later that the stability properties of each differ.
Note that $\phi= 0.5$ constitutes a line of discontinuities of the map, due to the sgn function that appears in the phase update rule. This discontinuity is clearly seen in the projections onto the $\phi-I$ plane (Fig. \ref{fig:period-phase-surfaces}C). The blue curve separates regions of $\phi-I$ space where the surface $z=F_2(I,\phi)$ lies above the diagonal plane $z=\phi$ or below it. In the projection  the region above the blue curve represents where the surface lies above the diagonal plane. Similarly, the region above the red curve is where the surface $z=F_1(I,\phi)$ lies below the diagonal plane $z=I$. 

\begin{figure}
\begin{center}
\includegraphics[width=0.45\textwidth]{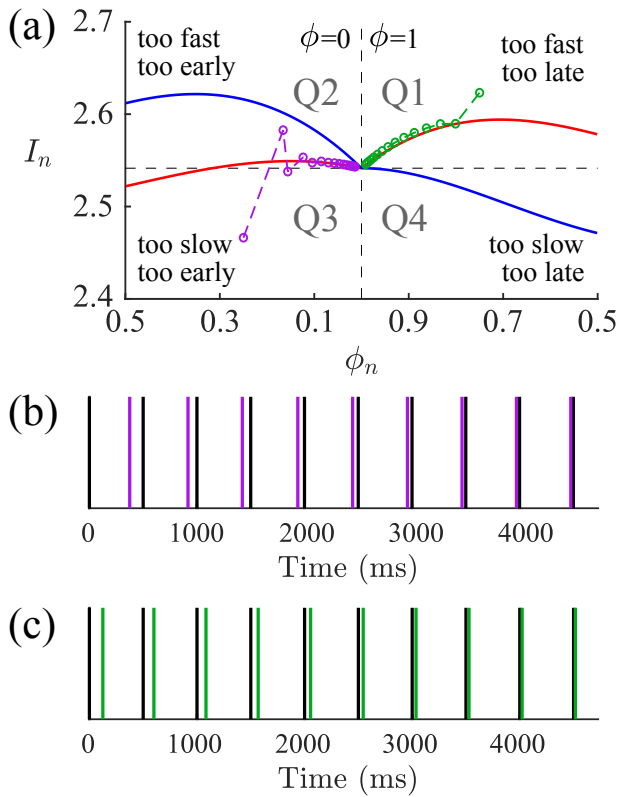}
\end{center}
\caption{\label{fig:order-preserving} \textbf{System can phase delay or phase advance during convergence to the synchronous solution.} (a) Recast phase plane with {\it dual} origin ($\phi=0$ and 1) corresponding to the synchronized solution. Note the non-standard ordering of the phase along the horizontal axis.  Iterates for two different initial conditions show that the phase is systematically decreased to $\phi=0$ when the BG is too early (purple) and increased toward $\phi=1$ when the BG is too late (green). 
(b) BG spike times (purple) and stimulus tone times (black) for the sequence of iterates that converges to $\phi=0$ via phase delay.
(c) BG spike times (green) and stimulus tone times (black) corresponding to the iterates that converges to $\phi=1$ by phase advance.
Parameter values given in Appendix \ref{appendix:parameters}.}
\end{figure}

It is useful to recast the $\phi-I$ phase plane by placing the origin at $I=I^*$ and at the {\it dual} values $\phi=0$ and 1 (Fig. \ref{fig:order-preserving}A), as both phase values correspond to the synchronized solution. Thus the entire vertical axis corresponds to both values $\phi=0$ and 1 and, as such, the BG spike and stimulus tone occurring at the same time. The horizontal axis is ordered to show how the phase changes, keeping in mind the dual role of the origin. In the left half-plane,  phase decreases (left to right) toward the origin.
in the right half-plane, phase increases (right to left) moving toward the origin.
This recast phase plane allows us to easily identify how $I$ and $\phi$ should be updated in each region, and in which direction the iterates should move. 
Namely, in the upper right  first quadrant (Q1), the phase indicates that the BG spiked after the stimulus tone, and is thus late. As a result the phase rule will try to speed up the BG.  However, the current value of $I_n$ is too large, thus the BG is too fast and will have an interspike interval less than $\tstim$. Thus the period rule will tend to slow down the BG. For example,  the iterates of the map with initial value $(\phi_0,I_0) = (0.75, 2.62)$ (green) systematically decrease the $I_n$ value to slow the BG down, while simultaneously increasing the phase towards the value 1 (i.e. convergence to the origin). Note that the iterates quickly converge to the red nullcline along which $I_{n+1}=I_n$. In a vicinity of this nullcline, the strengths of the phase and period rules are relatively balanced allowing the convergence towards synchrony to be monotone in phase.  The initial condition $(\phi_0,I_0) = (0.25, 2.47)$ (purple) lying in Q3 (lower left) corresponds to the BG being initially too slow, but also too early. Thus the period rule needs to increase the $I_n$ value, to speed up the BG, while the phase rule seeks to slow it down to decrease the phase. Again, the two-step updates eventually come into balance as can be seen the iterates move back and forth between Q2 and Q3 for a number of iterates, before converging towards the origin from Q2 (upper left) along the red nullcline. 

The time courses for the two examples clearly display how the system convergences towards the synchonized solution in an order preserving manner (Fig. \ref{fig:order-preserving}(a) and (c)). Initial condition lying in Q3 shows a systematic phase delay in its time course (Fig. \ref{fig:order-preserving}(b)); the BG spike time (purple) moves towards the stimulus tone time (black) that follow it. Note that due to the way phase is defined, a delay implies that the phase {\it decreases}. If the initial condition lies in Q1, then it is in the basin of attraction of the $\phi=1$ fixed point, and the BG spike times systematically phase advance (Fig. \ref{fig:order-preserving}(c)); the BG spike time (green) moves towards the stimulus tone time (black) that precedes it.

\subsection{Assessing stability of fixed points through linearization}
As demonstrated in Fig. \ref{fig:order-preserving}, depending on initial values, iterates of the map can converge to either the $\phi=0$ or $\phi=1$ fixed point. The convergence properties depend on the learning rule parameters $\delta_T$ and $\delta_\phi$ as well as $\tstim$. To better understand how the dynamics depend on these parameters, we assess the stability of the fixed points of the two-dimensional map by computing the eigenvalues of the Jacobian matrix,
\begin{equation}
J = \begin{pmatrix}
1+\delta_Tg(I^*) & -\delta_\phi(1-\delta_Tg(I^*) (1-\phi^*))  \\
-\frac{1}{\tstim}g(I^*) & 1+\frac{\delta_\phi}{\tstim}g(I^*)(1-\phi^*))
\end{pmatrix},
\end{equation}
where $(I^*,\phi^*)$ denotes a fixed point and $g(I)=dT(I)/dI$. 
The eigenvalues of $J$ are given by
\begin{align}
\label{eq:evalues}
\lambda_\pm &=  1+ \frac{1}{2}(\delta_T+\frac{\delta_\phi}{\tstim}(1-\phi^*))g(I^*) \\
&\pm\frac{1}{2} \sqrt{(\delta_T+\frac{\delta_\phi}{\tstim}(1-\phi^*))^2 g(I^*)^2+4\frac{\delta_\phi}{\tstim}g(I^*)}. \nonumber
\end{align}

\begin{figure}
\hspace{0.2em}
\begin{center}
\includegraphics[width=0.45\textwidth]{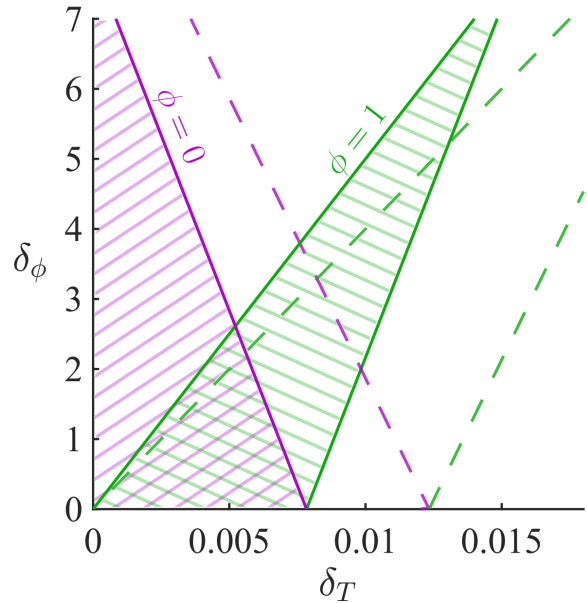}
\end{center}
\caption{\label{fig:stability-and-zero} \textbf{Fixed points have different stability properties.} The solid purple line shows the stability boundary for $\lambda=-1$ for the fixed point at $\phi^*=0$ with $\tstim=500$ ms. The fixed point is stable to the left of the purple curve. The solid green lines are stability boundaries for $\phi^*=1$ where the boundary emanating from the origin is for the complex eigenvalue condition $|\lambda_\pm|=1$ and the other curve is for $\lambda=-1$.  The fixed point is stable for parameter pairs lying between these two boundaries. The dashed line indicates how the boundaries change when $\tstim$ is reduced to $400$ ms.}
\end{figure}

A fixed point of the map is stable if both eigenvalues at the linearization lie inside the unit circle. We computed the stability boundaries for a fixed stimulus period $\tstim$ in the $\delta_T-\delta_\phi$ parameter plane (Fig. \ref{fig:stability-and-zero}) by solving $|\lambda|=1$ \eqref{eq:evalues}. This condition is met when an eigenvalue is real and equals $\pm 1$ or is complex with magnitude equal to 1. Although the fixed points at $\phi^*=0$ and $\phi^*=1$ correspond to the synchronous solution, they have different stability characteristics.  
The solid green lines correspond to the fixed point $\phi^*=1$ for $\tstim=500$. The line emanating from the origin is where the eigenvalues are complex with magnitude equal to one; the other when one of the eigenvalues equals -1. The solid purple line correspond to $\phi^*=0$ for $\tstim=500$ when an eigenvalue equals -1. The shaded green (purple) region correspond to parameter values for which the $\phi=1$ ($\phi=0$) fixed point is stable. It is straightforward to show that $dg/dI >0$ which implies that the stability boundary curves for $\lambda=-1$ shift to the right and for $|\lambda|=1$ have smaller slope as $\tstim$ is decreased. In turn, this implies that there is a larger set of parameters for which the $\phi=0$ and 1 fixed points are both stable when $\tstim$ is decreased (dashed lines in Fig. \ref{fig:stability-and-zero}).

Contrary to convention, stability of the fixed points is not sufficient to determine whether iterates of the map will actually converge to them. The reason is that the 2-dimensional map is built under the order preserving assumption, namely that every BG spike is followed by a stimulus tone. When one of the fixed points, say $\phi=1$, is a stable spiral point, the iterates of the map spiral in towards $\phi=1$, but when they cross the vertical axis there is a switch to convergence towards the $\phi=0$ fixed point. Such a switch corresponds to two consecutive BG spikes and hence, a violation of the 1-to-1 assumption of the map. In the next section we address how to handle these types of situations.

\section{\label{section-oiebmap}Order-indeterminant 2-dimensional map}

\begin{figure*}
\includegraphics[width=0.9\textwidth]{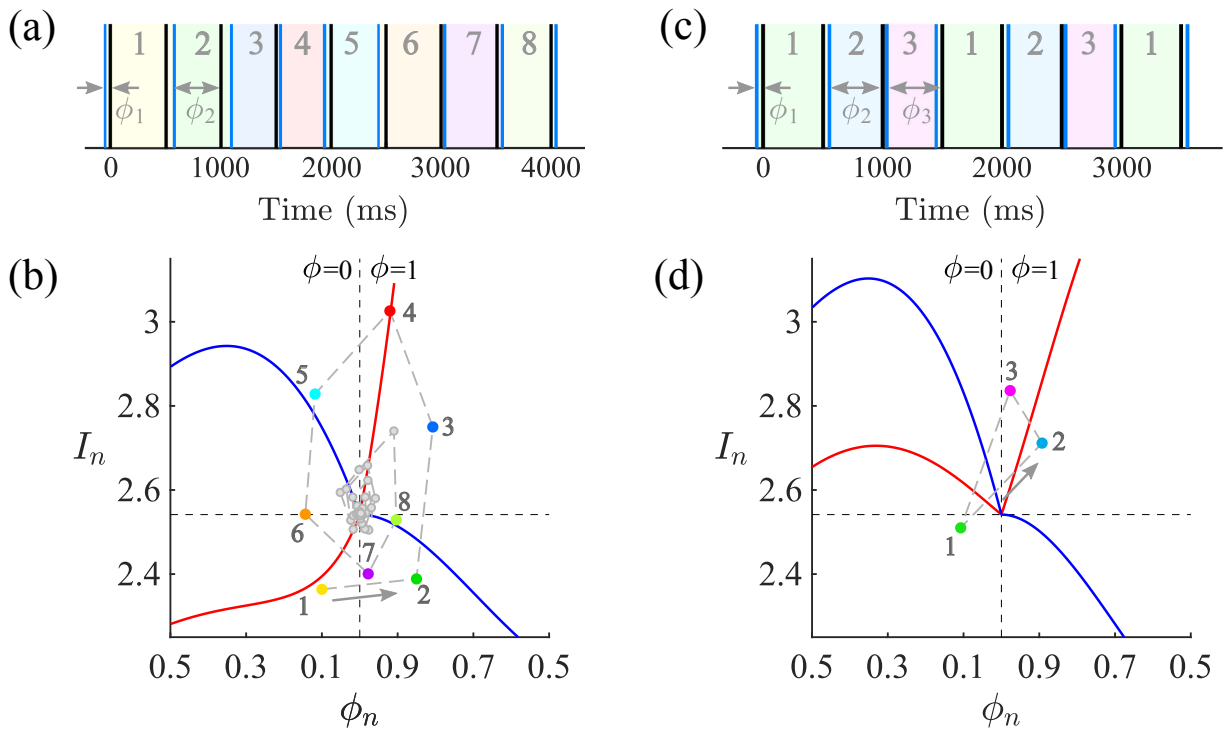}\\
\caption{\label{fig:order-indeterminant} \textbf{Order of BG spikes and stimulus tones is not always 1-to-1.} (a) Demonstration of order switching behaviour in the transient to synchronization. The BG switches between being to early and too late, before ultimately converging to the synchronous state where the BG spikes and stimulus tones occur at the same time. (b) $I-\phi$ phase plane for (a) depicting the phase and $I$ value at the start of each cycle. A order switch corresponds to the iterates crossing the $\phi=0$ -- $\phi=1$ line. Crossing from the $\phi=0$ side into the $\phi=1$ signifies two consecutive stimulus tones, while crossing from $\phi=1$ to $\phi=0$ indicates two consecutive BG spikes. 
(c) Time course for a period-3 solution. (d) Phase plane corresponding to (c).
Parameter values given in Appendix \ref{appendix:parameters}.}
\end{figure*}
As suggested in the previous section, convergence to the synchronous solution is not always order preserving. For example, if the BG spikes just before the stimulus tone, the phase correction may decrease $I$ too much and slow the BG down such that it does not fire again until after the next stimulus tone (Fig. \ref{fig:order-indeterminant}(a); cycle 1 and 6). Similarly, if the period rule over-corrects when speeding up the BG, the BG may fire again before the next stimulus tone  (Fig. \ref{fig:order-indeterminant}(a); cycle 4). 

As with the order preserving map, in the order-indeterminant case we say that a cycle starts and ends at a BG spike. However, the number of stimulus events per cycle can vary. For example, if there are two consecutive BG spikes, then there would be no stimulus spikes in this cycle, and if there are two consecutive stimulus tones, then both tones are said to be in the same cycle. The solution shown in Fig. \ref{fig:order-indeterminant}(c) and (d) is a period-3 solutions, where the first cycle (green) contains two stimulus events, the second (blue) contains 1 stimulus event and the third (pink) contains no stimulus events. 

At the beginning of each cycle, we compare the expected period of the BG (the period should no further updates happen) to the amount of time to the next stimulus spike. If the expected period is greater than this amount of time, the BG spike is followed by a stimulus tone, if not, it is followed by another BG spike. This leads us to define a function $H_S$ which will be zero if the expected period is less than the time to the next stimulus tone, and 1 otherwise,
\begin{equation}
H_S(I,\phi) = \Theta(T(I)-T_{stim}\phi),
\label{eq:BG-spike}
\end{equation}
where $\Theta (x)$ is the Heaviside function equal to 0 for $x<0$ and to 1 for $x\geq0$ and $T(I)$ is given by \eqref{eq:T}.
If there is a second stimulus tone in a cycle, we do not make another phase update, as we do not have sufficient information about the phase at this time. The BG has not spiked yet, so although we know that it is too late, we do not know by how much. Once the BG spikes again the period rule will act to speed up the BG. As such, there is at most one phase update per cycle. 

Using \eqref{eq:BG-spike} we write the order-indeterminant version of \eqref{eq:2d-map-i} as
\begin{align}
\label{eq:OIEB-i}
F_1(I,\phi) = I+ H_S(I,\phi)\delta_\phi \Delta I_\phi(\phi) + \delta_T\Delta I_T (I,\phi).
\end{align}
For the phase variable $\phi$, we note that if there are two or more stimulus spikes in a particular cycle, \eqref{eq:2d-map-phi} updates the phase to a negative number. We include a modulo 1 function to convert this to the proportion of time until the next stimulus tone rather than from the previous tone,
\begin{align}
\label{eq:OIEB-phi}
F_2(I,\phi) = \text{mod}\left(\phi+\frac{T_{stim}-T_n(I,\phi)}{T_{stim}},1\right).
\end{align}

With the OIEB map now defined, we return to describe the observed dynamics of Fig. \ref{fig:order-indeterminant}. During the transient to synchronization, the order of BG spikes and stimulus tones may alternate (Fig. \ref{fig:order-indeterminant}(a)). The BG is initially too slow and early, at the first stimulus spike ($t=0$ms) phase correction acts to slow down the BG further. At the next stimulus event ($t=500$ms), no updates occur as the phase has not been updated. Shortly afterwards, the BG spikes and $I$ is increased to speed up the BG and the phase is updated to $\phi\approx 0.8$. This marks the end of the first cycle (yellow). The second cycle contains 1 stimulus tone, and as $\phi\approx 0.8$, phase correction increases $I$, speeding up the BG. On the 4th cycle (red), the BG is too fast and the BG spikes again before the next stimulus event. Hence, there is no phase correction in this cycle. The order switches again in the 6th cycle (orange), when the BG is too slow and there are two consecutive stimulus tones. As with the first cycle, phase correction occurs at the first stimulus tone in the cycle, but not the second. The iterates of the map can be viewed on the $\phi-I$-plane (Fig. \ref{fig:order-indeterminant}(b)). In general, when an order switch occurs, the iterates cross from the left plane to the right plane for consecutive stimulus tones or from right to left for consecutive BG spikes. Thus the number of transitions between half-planes is equal to the number of order switching events. The iterates $I_n$ and $\phi_n$ correspond to the value of $I$ and $\phi$ at the start of the $n$th cycle, i.e. the value at the $n$th BG spike after period correction has been applied. The 8 cycles shown in the time course (Fig. \ref{fig:order-indeterminant}(a)) are colored accordingly in the phase plane and the remaining iterates before convergence to the synchronous solution are shown in grey.

The map also exhibits limit cycle behaviour (Fig. \ref{fig:order-indeterminant}(c)--(d)). A period-3 solutions is shown, where the first cycle of the periodic orbit (green) contains 2 stimulus tones, the second (blue) contains 1 stimulus tone and the third doesn't contain any stimulus tones (pink). At the first stimulus tone, the BG is too early and slightly too slow ($I_1<I^*$, $\phi_1<0.5$). The phase correction rule decreases $I$ and slows down the BG. As a result, the BG doesn't spike again until after the next stimulus tone, and there are  two stimulus tones in this cycle. The cycle ends at the next BG spike, after period correction is applied and $I$ is increased such that $I_2>I_1$. The BG is now too fast, but also late ($\phi > 0.9$). In the second cycle, phase correction increases the value of $I$ and period correction decreases it. We note that the phase correction is stronger as $I_3>I_2$. With a large phase $\phi \approx 1$ and $I_3>I^*$, the BG spikes again before the next stimulus tone. Hence, phase correction is not applied and the period rule acts to slow down the BG. The period-3 cycle then repeats.

\subsection{Dynamics of the OIEB map}
The values $(I^*(\tstim),0)$ and $(I^*(\tstim),1)$ are also fixed points of the OIEB map. Due to the discontinuity induced by the Heaviside function, it is not possible to linearize about the fixed points. Thus we cannot directly obtain stability information. 
However, the OIEB map reduces to the order preserving map of Section \ref{section-exactmap} whenever a sequence of BG spikes and stimulus tones occur consecutively. Hence, we can estimate the stability of the fixed points of the OIEB map using the order preserving map.

\begin{figure}
\hspace{0.5cm}
\begin{center}
\includegraphics[width=0.45\textwidth]{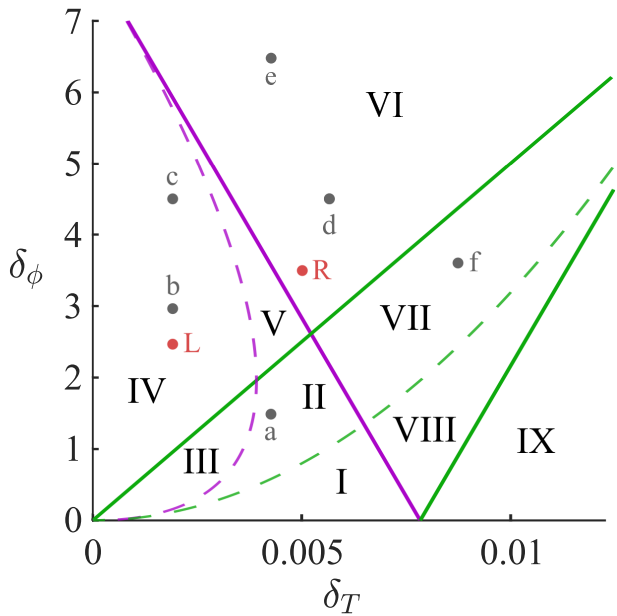}
\end{center}
\caption{\label{fig:deltaT-deltaP} \textbf{Stability characteristics of the order preserving map.} Solid lines indicate the stability boundaries, while dashed line show where the fixed points transition from nodes to spirals. The purple lines correspond to the fixed point at $\phi=0$, while the green are for the fixed point at $\phi=1$. The phase plane is divided into 9 distinct regions, whose stability characteristics are listed in Table \ref{tab1}. The red points correspond to the parameter value used in Fig. \ref{fig:order-indeterminant} (L -- left, R -- right), and grey points mark the parameter values used for the phase planes of the OIEB map shown in Fig. \ref{fig:gallery-behaviours}.}
\end{figure}

When the fixed point is a stable node, we see monotone convergence to $\phi=0$ or $\phi=1$, but when the fixed point is a spiral the iterates jump from converging towards $\phi=0$ to converging towards $\phi=1$ when there are two consecutive stimulus spikes and from converging to $\phi=1$ to converging to $\phi=0$ when there are consecutive BG spikes.
To quantify this we calculate where the discriminant of the eigenvalue given in (\ref{eq:evalues}) is 0 for $\tstim=500$ ms (Fig. \ref{fig:deltaT-deltaP}). The solid lines correspond to the stability boundaries, as in Fig. \ref{fig:stability-and-zero}, while the dashed curves represent the real-complex boundaries. For $\phi=0$ (purple) and $\phi=1$ (green), the eigenvalues are real to the right of the dashed curves. The curves separate the $\delta_T$-$\delta_\phi$ plane into 9 regions, labelled I-IX (Table \ref{tab1}). 

\begin{table}
\caption{Stability characteristics of fixed points of order preserving map \label{tab1}}
\begin{center}
\begin{tabular}{ccc}
\hline\noalign{\smallskip}
\multicolumn{1}{c}{Region} & \multicolumn{1}{c}{$\phi =0$} & 
\multicolumn{1}{c}{$\phi =1$}\\
\noalign{\smallskip}\hline\noalign{\smallskip}
\\[-.09in]
I & stable node & stable node\\
II & stable node & stable spiral\\
III & stable spiral & stable spiral\\
IV & stable spiral & unstable spiral\\
V & stable node & unstable spiral\\
VI & unstable node & unstable spiral\\
VII & unstable node & stable spiral\\
VIII & unstable node & stable node\\
IX & unstable node & unstable node\\
\noalign{\smallskip}\hline
\end{tabular}
\end{center}
\end{table}

\begin{figure*}
\begin{center}

\includegraphics[width=0.95\linewidth]{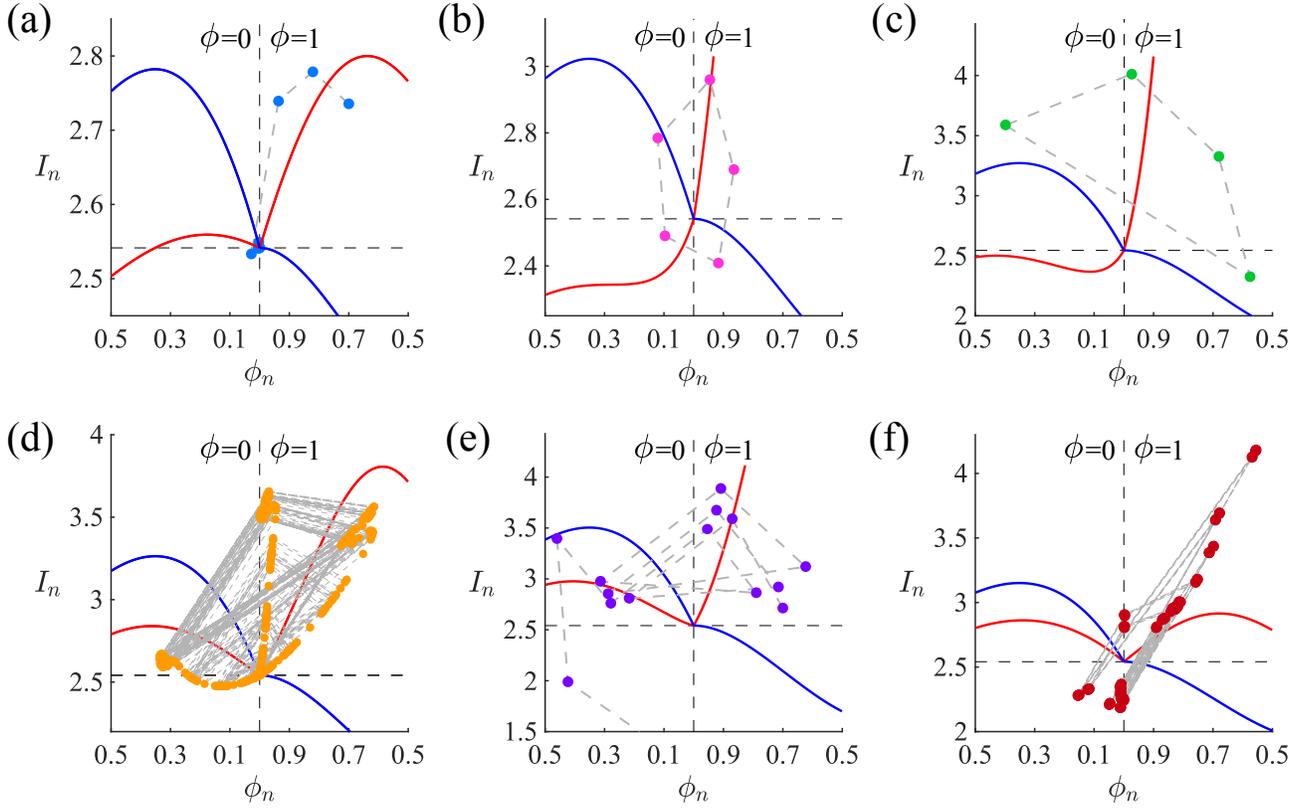}
\end{center}
\caption{\textbf{OIEB map exhibits a rich set of dynamics.} Phase planes for the points labelled a--f in Fig. \ref{fig:deltaT-deltaP}. (a) Iterates originally spiral in towards the fixed point at $\phi=1$, but there are two consecutive BG spikes and the iterates cross the $\phi=0$--$\phi=1$ line, and then converge to the synchronous solution from the left-half plane. (b) A period 5 solution where the iterates move counter-clockwise around the fixed point at $\phi=0$ and then $\phi=1$. (c) A period 4 solution where the number of BG spikes (4) and stimulus tones (3) are unequal.  (d) A chaotic attractor where $I_n$ and $\phi_n$ remain bounded but do not converge to either the synchronous solution or a periodic orbit. (e) Divergent behaviour, phase correction acts to slow down the BG but reduces $I_n$ below 1. (f) A period 104 orbit where order switching is frequent (many crossings of the $\phi=0$--$\phi=1$ line). Parameter values given in Appendix \ref{appendix:parameters}.}
\label{fig:gallery-behaviours}
\end{figure*}

Region I is the only part of parameter space where the eigenvalues of both fixed points are real and of magnitude less than one. As a result, it is the only region where convergence to a fixed point is consistently order-preserving, as seen in Fig. \ref{fig:order-preserving}.
In region II, initial condition in the basin of attraction of $\phi=0$ will exhibit order-preserving convergence, but those in the basin of attraction of $\phi=1$ will not. As $\phi_n$ increases through 1, it is reset to close to 0 (two consecutive BG spikes) and then converges monotonically to $\phi=0$ from there (Fig. \ref{fig:gallery-behaviours} (a)).
Both fixed points are stable spirals in region III, and iterates initially spiral in towards one of the two fixed points, but when two consecutive BG spikes or stimulus tones occur the iterates cross into the other plane and spiral towards the other fixed point as seen in Fig \ref{fig:order-indeterminant}(a)--(b). These crossings continue to occur until the system ultimately converges to the synchronized solution.

Moving into region IV and V, the fixed point at $\phi=1$ becomes unstable. In region IV, for parameter values suitably close to the stability boundary for $\phi=1$, the system behaviour is similar to that seen in region III, the iterates switch back and forth between the two fixed points but ultimately converge to the synchronized solution. If $\delta_\phi$ is increased further, we observe periodic orbits (Fig. \ref{fig:gallery-behaviours}(b)). This is a period-5 orbit with 2 order switches. The iterates orbit $\phi=0$ in a counter-clockwise direction (starting with iterate in Q2), until phase correction acts to slow down an already too slow BG, resulting in two consecutive stimulus event (iterate moves from left-half to right-half plane). The iterates switch to orbiting $\phi=1$, also in a counter-clockwise direction. Also in this region, periodic orbits with unequal numbers of BG spikes and stimulus tones exist (Fig. \ref{fig:gallery-behaviours}(c)). This solution corresponds to a period-4 orbit with 4 BG spikes, but only 3 stimulus tones. In general, a periodic orbit can have the number of BG spikes differ from the number of stimulus tones by at most one. In region V, the system exhibits behaviour akin to that of region II close to the $\phi=1$ stability boundary and limit cycle behaviour as we move away from the boundary.

Both fixed points are unstable in region VI. Close to the stability boundaries limit cycle behaviour exists. However, as $\delta_\phi$ is increased,  chaotic solutions arise (Fig. \ref{fig:gallery-behaviours}(d)). The iterates follow similar trajectories, but never return to the same point twice. In this region, further increasing $\delta_\phi$ leads to divergence, where $I$ is decreased below 1 and the BG stops oscillating (Fig. \ref{fig:gallery-behaviours}(e)).

In regions VII and VIII, the $\phi=1$ fixed point is stable, but the $\phi=0$ fixed point is not. For parameter choices close to the $\phi=0$ stability boundary the system converges to $\phi=1$ with a number of order switches in the transient. As we choose parameters away from the stability boundary, limit cycles emerge. In particular, we see large period limit cycles (Fig. \ref{fig:gallery-behaviours}(f)), with periods on the order of 100 or larger. Further increasing $\delta_T$ or $\delta_\phi$ in these regions once again  leads to chaotic attractors and then divergence when the minimum of $F_1(I,\phi)$ falls below 1. 
Both fixed points are unstable in region IX. As in region VI, we observe limit cycles close to the $\phi=1$ stability boundary. As $\delta_T$ is increased the limit cycles become chaotic attractors. Indeed, this follows naturally from our results on chaotic dynamics in the one-dimensional map of Section \ref{period-map}. Recall that when the fixed point of that map lost stability, there exists a set of $\delta_T$ values  over which period doubling bifurcations occur leading to chaos. For fixed $\tstim$ that set lies on the $\delta_{\phi}=0$ axis in Region IX of Fig. \ref{fig:deltaT-deltaP}. Thus for $\delta_\phi$ sufficiently small it is reasonable to intuit that such period doubling routes to chaos persist. Finally,  when $\delta_T$ is increased further the system diverges as minimum of $F_1(I,\phi)$ falls below 1. 

The above analysis demonstrates that the linearization at the fixed points of the order preserving map are useful in predicting and explaining the dynamics of the OIEB map. In particular, the characterization of the stability characteristics shown in Table \ref{tab1} allowed us to explain why certain iterates spiralled around the {\it dual} origin and either ended up converging as Fig. \ref{fig:order-indeterminant}(a) or remained bounded away from the origin as in Fig. \ref{fig:gallery-behaviours}(b). While we have not exhaustively explored parameter space, a more detailed parsing of parameter space would reveal additional boundaries that separate different types of behaviour within the same regions, e.g. a boundary curve separating the behaviour in Region V from chaotic \ref{fig:order-indeterminant}(d) to divergent \ref{fig:order-indeterminant}(e).

\section{Discussion}
In this paper, we have derived and analyzed a two-dimensional map that addresses how a beat generator (BG) neuronal network can learn the period and phase of an external periodic signal, such as a metronome. This map represents a type of error-correction algorithm \citep{Mates1994,Mates1994a,Repp2005,Repp2013} in which we define separate learning rules for period and phase that in turn adjust a biophysical parameter, $I$, of the BG to allow it to synchronize to the stimulus tones. The map originates from our prior work \citep{Bose2019} where we introduced an error correction algorithm that relied on the dynamics of neuronal oscillators. Here, we have simplified that presentation to create an event-based map, where iterations of the map correspond to cycles of the BG. 
Given that the order of BG spikes and stimulus events may alternate, we devised the map independent of the order of events. This led to what we termed an order-indeterminant event-based (OIEB) map. At each event, whether it be a BG spike or stimulus tone, the current values of $I$ and $\phi$ are used to determine what the next event is. Period correction is said to occur at BG spikes, while phase correction occurs at stimulus tones. As an iteration of the map is defined as a cycle of the BG, period correction occurs on every iteration, but phase correction may not.

Our analysis showed that there exists only a small region of $(\delta_T, \delta_\phi)$ (learning rate) parameter space (Region I in Fig. \ref{fig:deltaT-deltaP}) in which the synchronous solution was stable and approachable monotonically.  In this region, order preserving convergence signifies that the system is learning a beat by sequentially speeding up (phase advance) or slowing down (phase delay). In other regions (e.g. IV), the BG converges to the synchronous solution but does so in an order-indeterminant way, much in the manner a human participant would when learning a beat. As the strength of the period or phase rule becomes stronger, the BG loses the capability of convergence to the synchronous solution. This loss has two key implications: 1. there is a balance between how much either the phase or period rule can contribute to stable synchrony, and 2. the learning rates cannot be too large, and as such, convergence cannot be too fast.  Clearly the ability and rate at which humans learn to keep a beat improves with training.  Thus it is plausible that the learning rules proposed here may be subject to some sort of longer term plasticity that is shaped by experience.

In the one-dimensional period correction map as well as the two-dimensional period and phase correction map, we established the existence of both periodic and chaotic solutions. Chaos can arise in strictly continuous systems, such as the logistic map, but also occurs due to discontinuities in the definition of the map, as in switching systems. The discontinuity at $\phi=0.5$ constitutes a switching manifold. It would be of interest to see if the well-developed bifurcation theory \cite{DiBernardo2008} of such systems is applicable to OIEB maps.  Chaos also occurs in systems that exhibit border collision bifurcations in which a discontinuity of the map passes through a border on the domain of the map \citep{Nusse1995, Gardini2010}. In the order preserving map, the discontinuity of the map is fixed at $\phi=0.5$ independent of parameters. The OIEB map presents further discontinuities as a result of the Heaviside and modulo functions. These discontinuities appear to occur at interior points on the domain of the map. However, as can be noted from Fig. \ref{fig:gallery-behaviours}, the orientation of the nullclines changes with parameters. Perhaps such changes in orientation are analogous to border collisions.

An alternate strategy for keeping a beat involves direct periodic forcing  to entrain a set of mutually excitatory resonant oscillators \cite{Loehr2011, Large2015}.  These systems typically rely on weak coupling assumptions \cite{Hoppensteadt2012} and can be analyzed by classical methods of phase response curves \cite{Winfree2001}. Another recently introduced possibility for learning a beat is to employ a pulse-forced adaptable competition system consisting of units that ramp up at adjustable rates to a fixed threshold to produce specific time intervals\cite{Egger2019}. This approach follows the error-correction paradigm in that the ramping rate is subject to a learning rule. Our approach here is distinct from either of these strategies. We  propose that the BG is a type of  adaptive oscillatory system that is trying to develop an internal representation to an external dynamic source, the metronome.  Our approach may be considered as a temporal analogue in the case of timing of repetitive events to Kalman filtering in visual systems that corrects the internal presentation so its representation matches the external world \cite{Rao1997}. 

The map considered in this paper assumed that the event times were known and could be calculated exactly. However, in the context of humans estimating time, exact time information may not be available.
An alternate strategy that we proposed in previous work\citep{Bose2019} relies on subdividing time intervals into smaller reference intervals and counting the number of reference intervals between events; a form of discrete, integer-based estimation of timing. 
From a neuronal network point of view, this could be implemented by counting the number of oscillatory cycles of a fast spiking neuron or neural population between beat generator spikes and stimulus tones.
This method of discrete-time estimation could be included in the map with the introduction of an additional variable for each time interval to be estimated. These kinds of systems remain to be studied in future mathematical work. 

\appendix
\section{Parameter values for figures}
\label{appendix:parameters}
\begin{tabular}{|P{1.6cm}P{1.3cm}|P{1.2cm}|P{1.2cm}|P{1.2cm}|}
\hline  & & $\tstim$ & $\delta_T$ & $\delta_\phi$ \\ \hline 
Figure 1 &(a) & -- & 0.005 & --\\
& (b) & 500 & -- & --\\
& (c) & -- & 0.005 & --\\\hline
Figure 2 & & -- & - & --\\\hline
Figure 3 & & 500 & 0.005 & 0.5\\\hline
Figure 4 & & 500 & 0.005 & 0.5\\\hline
Figure 5 & \emph{solid} & 500 & -- & --\\
 & \emph{dashed} & 400 & -- & --\\
\hline
Figure 6 & (a)& 500 & 0.002 & 2.5\\
 & (b)& 500 & 0.002 & 2.5\\
 & (c) & 500 & 0.005 & 3.5\\
 & (d) & 500 & 0.005 & 3.5\\
\hline
Figure 7 &  & 500 & -- & --\\
\hline
Figure 8 & (a)& 500 & 0.0045 & 1.5\\
 & (b) & 500 & 0.002 & 3\\
 & (c) & 500 & 0.002 & 4.5\\  
 & (d) & 500 & 0.0055 & 4.5\\
 & (e) & 500 & 0.0045 & 6.5\\
 & (f) & 500 & 0.008 & 3.8\\
\hline
\end{tabular}

\newpage

\noindent \textbf{References} 
\bibliography{bib-file}

\begin{thebibliography}{27}%
\makeatletter
\providecommand \@ifxundefined [1]{%
 \@ifx{#1\undefined}
}%
\providecommand \@ifnum [1]{%
 \ifnum #1\expandafter \@firstoftwo
 \else \expandafter \@secondoftwo
 \fi
}%
\providecommand \@ifx [1]{%
 \ifx #1\expandafter \@firstoftwo
 \else \expandafter \@secondoftwo
 \fi
}%
\providecommand \natexlab [1]{#1}%
\providecommand \enquote  [1]{``#1''}%
\providecommand \bibnamefont  [1]{#1}%
\providecommand \bibfnamefont [1]{#1}%
\providecommand \citenamefont [1]{#1}%
\providecommand \href@noop [0]{\@secondoftwo}%
\providecommand \href [0]{\begingroup \@sanitize@url \@href}%
\providecommand \@href[1]{\@@startlink{#1}\@@href}%
\providecommand \@@href[1]{\endgroup#1\@@endlink}%
\providecommand \@sanitize@url [0]{\catcode `\\12\catcode `\$12\catcode
  `\&12\catcode `\#12\catcode `\^12\catcode `\_12\catcode `\%12\relax}%
\providecommand \@@startlink[1]{}%
\providecommand \@@endlink[0]{}%
\providecommand \url  [0]{\begingroup\@sanitize@url \@url }%
\providecommand \@url [1]{\endgroup\@href {#1}{\urlprefix }}%
\providecommand \urlprefix  [0]{URL }%
\providecommand \Eprint [0]{\href }%
\providecommand \doibase [0]{http://dx.doi.org/}%
\providecommand \selectlanguage [0]{\@gobble}%
\providecommand \bibinfo  [0]{\@secondoftwo}%
\providecommand \bibfield  [0]{\@secondoftwo}%
\providecommand \translation [1]{[#1]}%
\providecommand \BibitemOpen [0]{}%
\providecommand \bibitemStop [0]{}%
\providecommand \bibitemNoStop [0]{.\EOS\space}%
\providecommand \EOS [0]{\spacefactor3000\relax}%
\providecommand \BibitemShut  [1]{\csname bibitem#1\endcsname}%
\let\auto@bib@innerbib\@empty
\bibitem [{\citenamefont {Grahn}(2012)}]{Grahn2012}%
  \BibitemOpen
  \bibfield  {author} {\bibinfo {author} {\bibfnamefont {J.~A.}\ \bibnamefont
  {Grahn}},\ }\href@noop {} {\bibfield  {journal} {\bibinfo  {journal} {Topics
  in Cognitive Science}\ }\textbf {\bibinfo {volume} {4}},\ \bibinfo {pages}
  {585} (\bibinfo {year} {2012})}\BibitemShut {NoStop}%
\bibitem [{\citenamefont {Henry}, \citenamefont {Herrmann},\ and\ \citenamefont
  {Grahn}(2017)}]{Henry2017}%
  \BibitemOpen
  \bibfield  {author} {\bibinfo {author} {\bibfnamefont {M.~J.}\ \bibnamefont
  {Henry}}, \bibinfo {author} {\bibfnamefont {B.}~\bibnamefont {Herrmann}}, \
  and\ \bibinfo {author} {\bibfnamefont {J.~A.}\ \bibnamefont {Grahn}},\ }\href
  {\doibase 10.1371/journal.pone.0172454} {\bibinfo
  {journal} {PLoS ONE}\ }\textbf {\bibinfo {volume} {12}},\ \bibinfo {pages}
  {1} (\bibinfo {year} {2017})\BibitemShut {NoStop}%
\bibitem [{\citenamefont {Merchant}\ \emph {et~al.}(2015)\citenamefont
  {Merchant}, \citenamefont {Grahn}, \citenamefont {Trainor}, \citenamefont
  {Rohrmeier},\ and\ \citenamefont {Fitch}}]{Merchant2015}%
  \BibitemOpen
  \bibfield  {author} {\bibinfo {author} {\bibfnamefont {H.}~\bibnamefont
  {Merchant}}, \bibinfo {author} {\bibfnamefont {J.}~\bibnamefont {Grahn}},
  \bibinfo {author} {\bibfnamefont {L.}~\bibnamefont {Trainor}}, \bibinfo
  {author} {\bibfnamefont {M.}~\bibnamefont {Rohrmeier}}, \ and\ \bibinfo
  {author} {\bibfnamefont {W.~T.}\ \bibnamefont {Fitch}},\ }\href {\doibase
  10.1098/rstb.2014.0093} Philosophical Transactions of the Royal Society B: Biological Sciences \textbf {\bibinfo {volume} {370}},\ \bibinfo {pages} {2013} (\bibinfo {year}
  {2015})\BibitemShut {NoStop}%
\bibitem [{\citenamefont {Repp}(2005)}]{Repp2005}%
  \BibitemOpen
  \bibfield  {author} {\bibinfo {author} {\bibfnamefont {B.~H.}\ \bibnamefont
  {Repp}},\ }\href@noop {} {\bibfield  {journal} {\bibinfo  {journal}
  {Psychonomic Bulletin and Review}\ }\textbf {\bibinfo {volume} {12}},\
  \bibinfo {pages} {969} (\bibinfo {year} {2005})}\BibitemShut {NoStop}%
\bibitem [{\citenamefont {Repp}\ and\ \citenamefont {Su}(2013)}]{Repp2013}%
  \BibitemOpen
  \bibfield  {author} {\bibinfo {author} {\bibfnamefont {B.~H.}\ \bibnamefont
  {Repp}}\ and\ \bibinfo {author} {\bibfnamefont {Y.~H.}\ \bibnamefont {Su}},\
  }\href@noop {} {\bibfield  {journal} {\bibinfo  {journal} {Psychonomic
  Bulletin and Review}\ }\textbf {\bibinfo {volume} {20}},\ \bibinfo {pages}
  {403} (\bibinfo {year} {2013})}\BibitemShut {NoStop}%
\bibitem [{\citenamefont {Mates}(1994{\natexlab{a}})}]{Mates1994}%
  \BibitemOpen
  \bibfield  {author} {\bibinfo {author} {\bibfnamefont {J.}~\bibnamefont
  {Mates}},\ }\href@noop {} {\bibfield  {journal} {\bibinfo  {journal}
  {Biological Cybernetics}\ }\textbf {\bibinfo {volume} {70}},\ \bibinfo
  {pages} {463} (\bibinfo {year} {1994}{\natexlab{a}})}\BibitemShut {NoStop}%
\bibitem [{\citenamefont {Mates}(1994{\natexlab{b}})}]{Mates1994a}%
  \BibitemOpen
  \bibfield  {author} {\bibinfo {author} {\bibfnamefont {J.}~\bibnamefont
  {Mates}},\ }\href@noop {} {\bibfield  {journal} {\bibinfo  {journal}
  {Biological Cybernetics}\ }\textbf {\bibinfo {volume} {70}},\ \bibinfo
  {pages} {475} (\bibinfo {year} {1994}{\natexlab{b}})}\BibitemShut {NoStop}%
\bibitem [{\citenamefont {Achuthan}\ and\ \citenamefont
  {Canavier}(2009)}]{Achuthan2009}%
  \BibitemOpen
  \bibfield  {author} {\bibinfo {author} {\bibfnamefont {S.}~\bibnamefont
  {Achuthan}}\ and\ \bibinfo {author} {\bibfnamefont {C.~C.}\ \bibnamefont
  {Canavier}},\ }\href@noop {} {\bibfield  {journal} {\bibinfo  {journal}
  {Journal of Neuroscience}\ }\textbf {\bibinfo {volume} {29}},\ \bibinfo
  {pages} {5218} (\bibinfo {year} {2009})}\BibitemShut {NoStop}%
\bibitem [{\citenamefont {Oprisan}(2010)}]{Oprisan2010}%
  \BibitemOpen
  \bibfield  {author} {\bibinfo {author} {\bibfnamefont {S.~A.}\ \bibnamefont
  {Oprisan}},\ }\href@noop {} {\bibfield  {journal} {\bibinfo  {journal}
  {Journal of Theoretical Biology}\ }\textbf {\bibinfo {volume} {262}},\
  \bibinfo {pages} {232} (\bibinfo {year} {2010})}\BibitemShut {NoStop}%
\bibitem [{\citenamefont {Strumillo}\ and\ \citenamefont
  {Ruta}(2002)}]{Strumillo2002}%
  \BibitemOpen
  \bibfield  {author} {\bibinfo {author} {\bibfnamefont {P.}~\bibnamefont
  {Strumillo}}\ and\ \bibinfo {author} {\bibfnamefont {J.}~\bibnamefont
  {Ruta}},\ }\href {\doibase 10.1109/51.993195} IEEE Engineering in Medicine and Biology Magazine \textbf
  {\bibinfo {volume} {21}},\ \bibinfo {pages} {62} (\bibinfo {year}
  {2002})\BibitemShut {NoStop}%
\bibitem [{\citenamefont {Hamed}\ and\ \citenamefont
  {Grizzle}(2013)}]{Hamed2013}%
  \BibitemOpen
  \bibfield  {author} {\bibinfo {author} {\bibfnamefont {K.~A.}\ \bibnamefont
  {Hamed}}\ and\ \bibinfo {author} {\bibfnamefont {J.~W.}\ \bibnamefont
  {Grizzle}},\ }\href@noop {} {\bibfield  {journal} {\bibinfo  {journal} {IEEE
  Transactions on Robotics}\ }\textbf {\bibinfo {volume} {30}},\ \bibinfo
  {pages} {365} (\bibinfo {year} {2013})}\BibitemShut {NoStop}%
\bibitem [{\citenamefont {Bose}, \citenamefont {Byrne},\ and\ \citenamefont
  {Rinzel}(2019)}]{Bose2019}%
  \BibitemOpen
  \bibfield  {author} {\bibinfo {author} {\bibfnamefont {A.}~\bibnamefont
  {Bose}}, \bibinfo {author} {\bibfnamefont {{\'{A}}.}~\bibnamefont {Byrne}}, \
  and\ \bibinfo {author} {\bibfnamefont {J.}~\bibnamefont {Rinzel}},\ }\href
  {\doibase 10.1371/journal.pcbi.1006450}PLoS Computational Biology \textbf {\bibinfo {volume} {15}},\
  \bibinfo {pages} {e1006450} (\bibinfo {year} {2019})\BibitemShut {NoStop}%
\bibitem [{\citenamefont {Leonov}(1962)}]{Leonov1962}%
  \BibitemOpen
  \bibfield  {author} {\bibinfo {author} {\bibfnamefont {N.~N.}\ \bibnamefont
  {Leonov}},\ }\href@noop {} {\bibfield  {journal} {\bibinfo  {journal} {Dolk.
  Acad. Nauk SSSR}\ }\textbf {\bibinfo {volume} {143}},\ \bibinfo {pages}
  {1038} (\bibinfo {year} {1962})}\BibitemShut {NoStop}%
\bibitem [{\citenamefont {{Di Bernardo}}\ \emph {et~al.}(2008)\citenamefont
  {{Di Bernardo}}, \citenamefont {Budd}, \citenamefont {Champneys},
  \citenamefont {Kowalczyk}, \citenamefont {Nordmark}, \citenamefont {Tost},\
  and\ \citenamefont {Piiroinen}}]{DiBernardo2008}%
  \BibitemOpen
  \bibfield  {author} {\bibinfo {author} {\bibfnamefont {M.}~\bibnamefont {{Di
  Bernardo}}}, \bibinfo {author} {\bibfnamefont {C.~J.}\ \bibnamefont {Budd}},
  \bibinfo {author} {\bibfnamefont {A.~R.}\ \bibnamefont {Champneys}}, \bibinfo
  {author} {\bibfnamefont {P.}~\bibnamefont {Kowalczyk}}, \bibinfo {author}
  {\bibfnamefont {A.~B.}\ \bibnamefont {Nordmark}}, \bibinfo {author}
  {\bibfnamefont {G.~O.}\ \bibnamefont {Tost}}, \ and\ \bibinfo {author}
  {\bibfnamefont {P.~T.}\ \bibnamefont {Piiroinen}},\ }\href@noop {} {\bibfield
   {journal} {\bibinfo  {journal} {SIAM Review}\ }\textbf {\bibinfo {volume}
  {50}},\ \bibinfo {pages} {629} (\bibinfo {year} {2008})}\BibitemShut
  {NoStop}%
\bibitem [{\citenamefont {Glass}(1991)}]{Glass1991}%
  \BibitemOpen
  \bibfield  {author} {\bibinfo {author} {\bibfnamefont {L.}~\bibnamefont
  {Glass}},\ }\href@noop {} {\bibfield  {journal} {\bibinfo  {journal} {Chaos:
  An Interdisciplinary Journal of Nonlinear Science}\ }\textbf {\bibinfo
  {volume} {1}},\ \bibinfo {pages} {13} (\bibinfo {year} {1991})}\BibitemShut
  {NoStop}%
\bibitem [{\citenamefont {Diekman}\ and\ \citenamefont
  {Bose}(2016)}]{Diekman2016}%
  \BibitemOpen
  \bibfield  {author} {\bibinfo {author} {\bibfnamefont {C.}~\bibnamefont
  {Diekman}}\ and\ \bibinfo {author} {\bibfnamefont {A.}~\bibnamefont {Bose}},\
  }\href {\doibase 10.1177/0748730416662965} Journal of Biological Rhythms \textbf {\bibinfo {volume} {31}}
  (\bibinfo {year} {2016})\BibitemShut {NoStop}%
\bibitem [{\citenamefont {Booth}, \citenamefont {Xique},\ and\ \citenamefont
  {{Diniz Behn}}(2017)}]{Booth2017}%
  \BibitemOpen
  \bibfield  {author} {\bibinfo {author} {\bibfnamefont {V.}~\bibnamefont
  {Booth}}, \bibinfo {author} {\bibfnamefont {I.}~\bibnamefont {Xique}}, \ and\
  \bibinfo {author} {\bibfnamefont {C.~G.}\ \bibnamefont {{Diniz Behn}}},\
  }\href@noop {} {\bibfield  {journal} {\bibinfo  {journal} {SIAM Journal on
  Applied Dynamical Systems}\ }\textbf {\bibinfo {volume} {16}},\ \bibinfo
  {pages} {1089} (\bibinfo {year} {2017})}\BibitemShut {NoStop}%
\bibitem [{\citenamefont {May}(1976)}]{May1976}%
  \BibitemOpen
  \bibfield  {author} {\bibinfo {author} {\bibfnamefont {R.~M.}\ \bibnamefont
  {May}},\ }\href@noop {} {\bibfield  {journal} {\bibinfo  {journal} {Nature}\
  }\textbf {\bibinfo {volume} {261}},\ \bibinfo {pages} {459} (\bibinfo {year}
  {1976})}\BibitemShut {NoStop}%
\bibitem [{\citenamefont {Feigenbaum}(1978)}]{Feigenbaum1978}%
  \BibitemOpen
  \bibfield  {author} {\bibinfo {author} {\bibfnamefont {M.~J.}\ \bibnamefont
  {Feigenbaum}},\ }\href@noop {} {\bibfield  {journal} {\bibinfo  {journal}
  {Journal of Statistical Physics}\ }\textbf {\bibinfo {volume} {19}},\
  \bibinfo {pages} {25} (\bibinfo {year} {1978})}\BibitemShut {NoStop}%
\bibitem [{\citenamefont {Nusse}\ and\ \citenamefont
  {Yorke}(1995)}]{Nusse1995}%
  \BibitemOpen
  \bibfield  {author} {\bibinfo {author} {\bibfnamefont {H.~E.}\ \bibnamefont
  {Nusse}}\ and\ \bibinfo {author} {\bibfnamefont {J.~A.}\ \bibnamefont
  {Yorke}},\ }\href@noop {} {\bibfield  {journal} {\bibinfo  {journal}
  {International Journal of Bifurcation and Chaos}\ }\textbf {\bibinfo {volume}
  {5}},\ \bibinfo {pages} {189} (\bibinfo {year} {1995})}\BibitemShut {NoStop}%
\bibitem [{\citenamefont {Gardini}\ \emph {et~al.}(2010)\citenamefont
  {Gardini}, \citenamefont {Tramontana}, \citenamefont {Avrutin},\ and\
  \citenamefont {Schanz}}]{Gardini2010}%
  \BibitemOpen
  \bibfield  {author} {\bibinfo {author} {\bibfnamefont {L.}~\bibnamefont
  {Gardini}}, \bibinfo {author} {\bibfnamefont {F.}~\bibnamefont {Tramontana}},
  \bibinfo {author} {\bibfnamefont {V.}~\bibnamefont {Avrutin}}, \ and\
  \bibinfo {author} {\bibfnamefont {M.}~\bibnamefont {Schanz}},\ }\href@noop {}
  {\bibfield  {journal} {\bibinfo  {journal} {International Journal of
  Bifurcation and Chaos}\ }\textbf {\bibinfo {volume} {20}},\ \bibinfo {pages}
  {3085} (\bibinfo {year} {2010})}\BibitemShut {NoStop}%
\bibitem [{\citenamefont {Loehr}, \citenamefont {Large},\ and\ \citenamefont
  {Palmer}(2011)}]{Loehr2011}%
  \BibitemOpen
  \bibfield  {author} {\bibinfo {author} {\bibfnamefont {J.~D.}\ \bibnamefont
  {Loehr}}, \bibinfo {author} {\bibfnamefont {E.~W.}\ \bibnamefont {Large}}, \
  and\ \bibinfo {author} {\bibfnamefont {C.}~\bibnamefont {Palmer}},\
  }\href@noop {} {\bibfield  {journal} {\bibinfo  {journal} {Journal of
  Experimental Psychology: Human Perception and Performance}\ }\textbf
  {\bibinfo {volume} {37}},\ \bibinfo {pages} {1292} (\bibinfo {year}
  {2011})}\BibitemShut {NoStop}%
\bibitem [{\citenamefont {Large}, \citenamefont {Herrera},\ and\ \citenamefont
  {Velasco}(2015)}]{Large2015}%
  \BibitemOpen
  \bibfield  {author} {\bibinfo {author} {\bibfnamefont {E.~W.}\ \bibnamefont
  {Large}}, \bibinfo {author} {\bibfnamefont {J.~A.}\ \bibnamefont {Herrera}},
  \ and\ \bibinfo {author} {\bibfnamefont {M.~J.}\ \bibnamefont {Velasco}},\
  }\href@noop {} {\bibfield  {journal} {\bibinfo  {journal} {Frontiers in
  Systems Neuroscience}\ }\textbf {\bibinfo {volume} {9}},\ \bibinfo {pages}
  {1} (\bibinfo {year} {2015})}\BibitemShut {NoStop}%
\bibitem [{\citenamefont {Hoppensteadt}\ and\ \citenamefont
  {Izhikevich}(2012)}]{Hoppensteadt2012}%
  \BibitemOpen
  \bibfield  {author} {\bibinfo {author} {\bibfnamefont {F.~C.}\ \bibnamefont
  {Hoppensteadt}}\ and\ \bibinfo {author} {\bibfnamefont {E.~M.}\ \bibnamefont
  {Izhikevich}},\ }\href@noop {} {\emph {\bibinfo {title} {{Weakly connected
  neural networks}}}},\ Vol.\ \bibinfo {volume} {126}\ (\bibinfo  {publisher}
  {Springer Science and Business Media},\ \bibinfo {year} {2012})\BibitemShut
  {NoStop}%
\bibitem [{\citenamefont {Winfree}(2001)}]{Winfree2001}%
  \BibitemOpen
  \bibfield  {author} {\bibinfo {author} {\bibfnamefont {A.~T.}\ \bibnamefont
  {Winfree}},\ }\href@noop {} {\emph {\bibinfo {title} {{The geometry of
  biological time}}}},\ Vol.~\bibinfo {volume} {12}\ (\bibinfo  {publisher}
  {Springer Science and Business Media},\ \bibinfo {year} {2001})\BibitemShut
  {NoStop}%
\bibitem [{\citenamefont {Egger}, \citenamefont {Le},\ and\ \citenamefont
  {Jazayeri}(2019)}]{Egger2019}%
  \BibitemOpen
  \bibfield  {author} {\bibinfo {author} {\bibfnamefont {S.~W.}\ \bibnamefont
  {Egger}}, \bibinfo {author} {\bibfnamefont {N.~M.}\ \bibnamefont {Le}}, \
  and\ \bibinfo {author} {\bibfnamefont {M.}~\bibnamefont {Jazayeri}},\ }\href
  {\doibase 10.1101/712141} bioRxiv\ ,\ \bibinfo {pages} {712141} (\bibinfo {year}
  {2019})\BibitemShut {NoStop}%
\bibitem [{\citenamefont {Rao}\ and\ \citenamefont {Ballard}(1997)}]{Rao1997}%
  \BibitemOpen
  \bibfield  {author} {\bibinfo {author} {\bibfnamefont {R.~P.~N.}\
  \bibnamefont {Rao}}\ and\ \bibinfo {author} {\bibfnamefont {D.~H.}\
  \bibnamefont {Ballard}},\ }\href@noop {} {\bibfield  {journal} {\bibinfo
  {journal} {Neural Computation}\ }\textbf {\bibinfo {volume} {9}},\ \bibinfo
  {pages} {721} (\bibinfo {year} {1997})}\BibitemShut {NoStop}%
\end{thebibliography}%

\end{document}